\documentclass[a4paper,11pt]{amsart}

\usepackage{amsfonts,amssymb,mathrsfs,amsmath,amscd,epsfig,setspace,amstext,color,graphics,verbatim}
\usepackage{paralist,hyperref}
\usepackage{frenchineq}

\usepackage{fullpage}

\usepackage{amsthm}

\usepackage{mathtools}
\usepackage{mathrsfs}
\usepackage{fixmath}
\usepackage{comment}

\usepackage{phdalgo}
\usepackage[capitalize]{cleveref}

\crefformat{equation}{\textup{#2(#1)#3}}
\crefrangeformat{equation}{\textup{#3(#1)#4--#5(#2)#6}}
\crefmultiformat{equation}{\textup{#2(#1)#3}}{ and \textup{#2(#1)#3}}
{, \textup{#2(#1)#3}}{, and \textup{#2(#1)#3}}
\crefrangemultiformat{equation}{\textup{#3(#1)#4--#5(#2)#6}}%
{ and \textup{#3(#1)#4--#5(#2)#6}}{, \textup{#3(#1)#4--#5(#2)#6}}{, and \textup{#3(#1)#4--#5(#2)#6}}

\Crefformat{equation}{#2Equation~\textup{(#1)}#3}
\Crefrangeformat{equation}{Equations~\textup{#3(#1)#4--#5(#2)#6}}
\Crefmultiformat{equation}{Equations~\textup{#2(#1)#3}}{ and \textup{#2(#1)#3}}
{, \textup{#2(#1)#3}}{, and \textup{#2(#1)#3}}
\Crefrangemultiformat{equation}{Equations~\textup{#3(#1)#4--#5(#2)#6}}%
{ and \textup{#3(#1)#4--#5(#2)#6}}{, \textup{#3(#1)#4--#5(#2)#6}}{, and \textup{#3(#1)#4--#5(#2)#6}}

\DeclareMathAlphabet{\mathbfcal}{OMS}{cmsy}{b}{n}
\DeclareMathAlphabet{\mathbbold}{U}{bbold}{m}{n}

\newcommand{\shapley}{F}
\newcommand{\bias}{v}

\providecommand{\arxiv}[2][]{\href{http://www.arxiv.org/abs/#2}{arXiv:#2}}

\DeclareMathOperator*{\val}{\mathsf{val}}

\newcommand{\ie}{i.e.}
\newcommand{\trop}[1][]{\ifthenelse{\equal{#1}{}}{ \mathbb{T}_{\max}}{ \mathbb{T}(#1) }}

\newcommand{\minstates}{n}
\newcommand{\maxstates}{m}
\newcommand{\chancestates}{q}
\newcommand{\lcw}{\underline{\mathrm{cw}}}
\newcommand{\ucw}{\overline{\mathrm{cw}}}
\newcommand{\cond}{\operatorname{cond}}

\newcommand{\R}{\mathbb{R}}

\newcommand{\N}{\mathbb{N}}

\newcommand{\tplus}{\oplus}  

\newcommand{\tdot}{\odot}

\newcommand{\strictleq}{\ll}

\newcommand{\zero}{{-\infty}}

\newcommand{\sC}{\mathbfcal{C}}
\newcommand{\sS}{\mathcal{S}}
\newcommand{\PP}{\mathscr{P}}
\newcommand{\mpzero}{\mathbbb{0}}

\newcommand{\rminb}{\overline{\mathbb{T}}_{\min}}
\newcommand{\defi}{\coloneqq}
\newcommand{\mybot}{\mathbf{b}}
\newcommand{\mytop}{\mathbf{t}}
\newcommand{\oneto}[1]{[#1]}

\newcommand{\HNorm}[1]{{\| #1 \|}_{\rm{H}}}

\newcommand{\puiseux}{\mathbb{K}}

\newcommand{\nnpuiseux}{\mathbb{K}_{\geq 0}}
\DeclareMathAlphabet{\mathbbold}{U}{bbold}{m}{n}
\newcommand{\bo}[1]{\mathbold{#1}}

\newcommand{\states}{n}
\newcommand{\state}{i}

\newcommand{\stateII}{l}
\newcommand{\comd}{M}
\newcommand{\potential}{Q}
\newcommand{\transition}{P}
\newcommand{\rstates}{k}
\newcommand{\nstates}{n}

\newcommand{\hilb}[1]{\| #1 \|_{\mathrm{H}}}

\newcommand{\Payoff}{W}
\newcommand{\invmes}{\pi}

\newtheorem{theorem}{Theorem}
\newtheorem{proposition}[theorem]{Proposition}
\newtheorem{corollary}[theorem]{Corollary}
\newtheorem{lemma}[theorem]{Lemma}
\theoremstyle{definition}
\newtheorem{assumption}{Assumption}
\newtheorem{definition}{Definition}
\theoremstyle{remark}
\newtheorem{remark}[theorem]{Remark}
\newtheorem{example}[theorem]{Example}

\usepackage[colorinlistoftodos,bordercolor=orange,backgroundcolor=orange!20,linecolor=orange,textsize=scriptsize]{todonotes}

\DeclareMathAlphabet{\mathbbb}{U}{bbold}{m}{n}

\begin{document}

\title{Condition numbers of stochastic mean payoff games and what they say about nonarchimedean semidefinite programming}
\date{\today}

\thanks{X. Allamigeon, S. Gaubert, and M. Skomra were partially supported by the ANR projects CAFEIN (ANR-12-INSE-0007) and MALTHY (ANR-13-INSE-0003), by the PGMO program of EDF and Fondation Math\'ematique Jacques Hadamard, and by
the ``Investissement d'avenir'', r{\'e}f{\'e}rence ANR-11-LABX-0056-LMH,
LabEx LMH. M.~Skomra is  supported by a grant from R{\'e}gion Ile-de-France.
The first two authors also gratefully acknowledge the support of Mittag-Leffler Institute.
}%

\author{Xavier Allamigeon}
\author{St{\'e}phane Gaubert}
\author{Ricardo D. Katz}
\author{Mateusz Skomra}

\address{\textit{X.~Allamigeon, S.~Gaubert, and M.~Skomra}: INRIA and CMAP, \'Ecole Polytechnique, CNRS, 91128 Palaiseau Cedex France}
\email{firstname.lastname@inria.fr}

\address{\textit{R.~D.~Katz}: CONICET-CIFASIS, Bv. 27 de Febrero 210 bis, 2000 Rosario, Argentina}
\email{katz@cifasis-conicet.gov.ar}

\keywords{Semidefinite programming, stochastic games, tropical geometry, nonarchimedean fields}

\begin{abstract}
Semidefinite programming can be considered over any real closed field,
including fields of Puiseux series equipped with their nonarchimedean
valuation. Nonarchimedean semidefinite programs encode parametric
families of classical semidefinite programs, for sufficiently
large values of the parameter. Recently, a correspondence has
been established between nonarchimedean semidefinite programs
and stochastic mean payoff games with perfect information. 
This correspondence relies on tropical geometry.
It allows one to solve generic nonarchimedean semidefinite feasibility problems, of large scale, by means of stochastic game algorithms. 
In this paper, we show that the mean payoff of these games can be interpreted
as a condition number for the corresponding nonarchimedean feasibility problems.
This number measures how close a feasible instance is from being infeasible,
and vice versa. We show that it coincides with the maximal radius of a ball in Hilbert's projective metric, that is included in the feasible set. The geometric interpretation of the condition number relies in particular on a duality theorem for tropical semidefinite feasibility programs. 
Then, we bound the complexity of the feasibility problem
in terms of the condition number. We finally give explicit bounds for
this condition number, in terms of the characteristics of the 
stochastic game. As a consequence, we show that the simplest
algorithm to decide whether a stochastic mean payoff game is winning,
namely value iteration, has a pseudopolynomial complexity
when the number of random positions is fixed. 
\end{abstract}

\maketitle

\section{Introduction}
\subsection{Motivation}
Semidefinite programming (SDP) consists in optimizing
a linear function over a spectrahedron, the latter
being the intersection of a cone of positive semidefinite
matrices with an affine space. Semidefinite programs arise
in a number of applications from engineering sciences
and combinatorial optimization. 
We refer the reader to~\cite{siam_parrilo,matousekbook} for more background on the theory
and applications of semidefinite programming.

Spectrahedra form a class of convex semialgebraic sets. Even though these sets are usually defined over the field of real numbers, their definition is meaningful over any real closed field. In particular, the complexity of SDP and related questions can be investigated over real closed nonarchimedean fields, like
fields of Puiseux series. Such nonarchimedean SDP instances, which arise in perturbation theory, encode 
parametric families of classical SDP instances (over the reals),
for large enough (or small enough)
values of the parameter. The study of the nonarchimedean 
case is also motivated by unsettled questions concerning the complexity of ordinary SDP. Indeed, the latter are 
solvable in ``polynomial time'' only in a restricted sense.
More precisely, complexity bounds for SDP, obtained by the ellipsoid or interior
point methods, are only polynomial
in the log of certain metric estimates whose bit-size can be doubly
exponential in the size of the input~\cite{deklerk_vallentin}. 
It is unknown whether the SDP feasibility problem
belongs to $\text{NP}$. 

Semidefinite feasibility problems over the nonarchimedean valued field of Puiseux series have been studied in~\cite{issac2016jsc}. It is shown there that, under a genericity condition, these problems are equivalent to solving stochastic mean payoff games with perfect information and finite state and action spaces.
Stochastic mean payoff games
have an unsettled complexity: they belong to $\text{NP} \cap \text{coNP}$ but
no polynomial time algorithm is currently known~\cite{condon,andersson_miltersen}. However, several
practically efficient algorithms to solve stochastic mean payoff
games have been developed. In this way, one can solve nonarchimedean
semidefinite instances of a scale probably unreachable by interior
point methods. For instance, the benchmarks presented in~\cite{issac2016jsc} show that random instances of these problems with as many as $10000$ variables
could be solved by value iteration in a few seconds. Hard instances
are experimentally concentrated in a small ``phase transition''
region of the parameter space.

\subsection{Main results}
In order to explain why value iteration
is so efficient on many nonarchimedean SDP instances, we introduce here a notion of {\em condition
number} for stochastic mean payoff games. Essentially,
for a feasible instance, the condition number is the inverse of the distance
of the data to an infeasible instance, and vice versa. We show
that this condition number coincides with the absolute
value of the mean payoff. 
We establish a universal bound for the time of convergence
of value iteration, involving 
the condition number and an auxiliary metric
estimate, the distance of point $0$ to the set
of ``bias vectors'' (\cref{Th:NitsBound}). Then, we effectively bound
the condition number and the latter distance, 
for stochastic mean payoff games with perfect information (\cref{th-cond,th-bias}). We arrive, in particular,
at a bound that becomes pseudopolynomial
when the number of ``random'' positions 
of the game is fixed.  

To arrive at these results, we develop
a metric geometry approach of the condition
number. We use Hilbert's
projective metric, which arises in
Perron--Frobenius theory~\cite{nussbaum88}.
The same metric, up to a logarithmic change
of variable, arises in tropical geometry~\cite{cgq02,polyhedra_equiv_mean_payoff}. We also prove duality results for stochastic mean payoff games,
showing, essentially, that the condition number
of the primal and dual problems coincide. 
In summary, our main results show that the complexity of
value iteration is governed by metric geometry
properties: this leads to a general method
to derive complexity bounds, 
which can be applied to various classes
of Shapley operators. 

\subsection{Related works}
When specialized to stochastic
mean payoff games with perfect information, our bounds should be compared with
the one of Boros, Elbassioni, Gurvich, and Makino~\cite{boros_gurvich_makino}.
The authors of~\cite{boros_gurvich_makino} generalize the 
``pumping'' algorithm, developed for deterministic games by
Gurvich, Karzanov, and Khachiyan~\cite{gurvich}, to the case of stochastic games.
The resulting algorithm is
also pseudopolynomial
if the number of random positions is fixed.
The algorithm of Ibsen-Jensen and Miltersen~\cite{ibsen-jensen_miltersen}
yields a stronger bound in the case of simple stochastic games,
still assuming that the number of random positions
is fixed. 

The duality results in \cref{sec-duality} 
extend to stochastic games some  duality results
for deterministic games by Grigoriev and Podolskii~\cite{GrigorievPodolskii}.
In contrast, our approach builds on~\cite{polyhedra_equiv_mean_payoff},
deriving duality results from a minimax Collatz--Wielandt type
theorem of Nussbaum~\cite{nussbaum}.
Other duality results, by Bodirsky and Mamino, in the context
of satisfiability problems, have appeared in~\cite{bodirsky_mamino}.

\subsection{Organization of the paper}
Earlier results on the relation between nonarchimedean semidefinite programming
and stochastic mean payoff games are presented in \cref{sec-motivation}, leading  to the introduction of the notion of condition number.
Some background on nonlinear Perron--Frobenius
theory is presented in \cref{sec:perron_frobenius}.
The new results are included in \cref{sec-metric}, in which
we characterize the condition number, and in \cref{sec-complexity}, in
which we derive complexity estimates for value iteration
in terms of the condition number. This is a preliminary announcement of the results. The proofs will appear in a subsequent version.

\section{Motivation: the correspondence between nonarchimedean semidefinite programming and stochastic mean payoff games}\label{sec-motivation} 
In this section, we summarize some of the main results of~\cite{issac2016jsc},
which motivate the present work. Throughout this paper, given $k\in \N$, we denote the set $\{1, \dots, k\}$ by~$\oneto{k}$.

\subsection{Nonarchimedean semidefinite programs}

We start by introducing semidefinite programming over nonarchimedean fields. More specifically, the model of nonarchimedean field used in this paper is the field $\puiseux$ of \emph{(absolutely convergent generalized real) Puiseux series}, which are series in the parameter $t$ of the form
\begin{equation}
\bo x = \sum_{i = 1}^{\infty} c_{\lambda_{i}}t^{\lambda_{i}} \, , \label{eq:series}
\end{equation}
where \begin{inparaenum}[(i)]
\item $(\lambda_{i})_{i \ge 1}$ is a strictly decreasing sequence of real numbers that is either finite or unbounded, 
\item $c_{\lambda_{i}} \in \R \setminus \{ 0\}$ for all $\lambda_{i}$,
\item and the series~\cref{eq:series} is absolutely convergent for $t \in \R$ sufficiently large. 
\end{inparaenum} There is also a special, empty series, which is denoted by $\mathbf 0$. 
The field $\puiseux$ is ordered, with a total order defined by $\bo x > \mathbf{0} \iff c_{\lambda_{1}} > 0$. In addition, it is known that $\puiseux$ is a real closed field~\cite{van_dries_power_series}. Actually, our approach applies
to other nonarchimedean fields with a real value group~\cite{tropical_spectrahedra}, but it is helpful to have a concrete field in mind, like
$\puiseux$.
Henceforth, we denote by $\nnpuiseux$ the set of nonnegative series, i.e., the series $\bo x$ that satisfy $\bo x = \mathbf{0}$ or $\bo x > \mathbf{0}$.

Given symmetric matrices $\bo Q^{(0)}, \bo Q^{(1)}, \dots, \bo Q^{(n)} \in \puiseux^{m \times m}$, we define the associated \emph{spectrahedron} (over Puiseux series) as the set
\begin{equation}\label{puiseux_spectrahedron}
\bigl\{ \bo x \in \puiseux^n \colon \bo Q^{(0)} + \bo x_1 \bo Q^{(1)} + \dots + \bo x_n \bo Q^{(n)} \; \text{is PSD} \bigr\} \, ,
\end{equation}
where ``PSD'' stands for \emph{positive semidefinite}. (We point out that the definition of positive semidefinite matrices makes sense over any real closed field.) The problem which we are interested in is to determine whether a spectrahedron over Puiseux series is empty or not. This corresponds to the analog of the semidefinite feasibility problem over the field $\puiseux$. This problem is also related to the standard semidefinite feasibility problem over the field of real numbers associated with the spectrahedra
\begin{equation}\label{real_puiseux_spectrahedron}
\bigl\{ x \in \R^n \colon \bo Q^{(0)}(t) + x_1 \bo Q^{(1)}(t) + \dots + x_n \bo Q^{(n)}(t) \; \text{is PSD} \bigr\} \, ,
\end{equation}
for $t$
large enough.
Here, $\bo Q^{(i)}(t)$ stands for the real symmetric matrix obtained by evaluating the entries of $\bo Q^{(i)}$ at the value $t$. The relation between the problem over Puiseux series and the one over real numbers is described in the following lemma, and is a consequence of quantifier elimination over real closed fields:
\begin{lemma}
The spectrahedron \cref{puiseux_spectrahedron} over the field $\puiseux$ is empty if and only if, for $t$ sufficiently large, the spectrahedron \cref{real_puiseux_spectrahedron} over $\R$ is empty.
\end{lemma}
In this paper, we consider a slightly different problem which already retains much of the difficulty of the semidefinite feasibility problem over the field $\puiseux$: given symmetric matrices $\bo Q^{(1)}, \dots, \bo Q^{(n)} \in \puiseux^{m \times m}$, determine whether the following \emph{spectrahedral cone}
\begin{equation}\label{puiseux_spectra_cone}
\bigl\{ \bo x \in \nnpuiseux^n \colon  \bo x_1 \bo Q^{(1)} + \dots + \bo x_n \bo Q^{(n)} \; \text{is PSD} \bigr\} \, 
\end{equation}
is \emph{trivial}, meaning that it is reduced to the zero point. We refer to~\cite{issac2016jsc} for further details on the relation between the original semidefinite feasibility problems and the problems above for spectrahedral cones.

\subsection{Valuation map and tropical semifield}\label{subsec:tropical}

As a nonarchimedean field, $\puiseux$ is equipped with a \emph{valuation} map $\val \colon \puiseux \to \R \cup \{ -\infty\}$ defined by $\val(\bo x) \defi \lambda_{1}$ for $\bo x \neq \mathbf{0}$ as in~\cref{eq:series}, and $\val(\mathbf{0}) \defi -\infty$. This valuation map has the following properties:
\begin{align}
\val(\bo x + \bo y) &\leq \max( \val (\bo x ), \val (\bo y))
\label{e-val}\\
\qquad 
\val(\bo x  \bo y) &= \val (\bo x ) + \val (\bo y) \, .
\label{e-val2}
\end{align}
We point out that equality holds in~\cref{e-val} as soon as the leading terms of $\bo x$ and $\bo y$ do not cancel. In particular, this condition is satisfied when $\bo x, \bo y\geq \mathbf{0}$.

The \emph{tropical} (or \emph{max-plus}) semifield
$\trop$ can be though of as the image of $\nnpuiseux$ by the valuation map. More precisely, this semifield is defined as the set $\trop \defi \R \cup \{-\infty \}$ endowed with the addition
$x \tplus y \coloneqq \max(x,y)$ and the multiplication $x \tdot y \coloneqq x + y$. The term ``semifield'' refers to the fact that the addition does not have an opposite law.  The reader may consult~\cite{bcoq,butkovic,maclagan_sturmfels} for more information on the tropical semifield.

The operations above are extended in the usual way to matrices with entries in $\trop$. The resulting matrix product is also denoted by $\odot$. 
Henceforth, for any $z\in \trop^n$ and $\beta \in \trop$, we denote by $\beta + z$ the vector of $\trop^n$ with entries $\beta + z_i$. Finally, we denote by $\mpzero$ the neutral element for addition in $\trop$ (i.e., $\mpzero \defi -\infty$), as well as any vector that has all components equal to $\mpzero$. 

We consider $\trop$ equipped with the topology defined by the distance $(a,b) \mapsto |\exp(a) - \exp (b)|$, and $\trop^n$ equipped with the product topology. On $\R^n$ we also use {\em Hilbert's seminorm}~\cite{cgq02}, defined by $\HNorm{x} \defi \mytop(x) - \mybot(x)$, where $\mytop (x) \defi \max_{i\in \oneto{n}} x_i$ and $\mybot (x) \defi \min_{i\in \oneto{n}} x_i$. This seminorm induces a norm on the quotient space of $\R^n$ by the tropical parallelism relation, which is defined by: $x \parallel y$ if, and only if, there exists $\alpha \in \R$ such that $x = \alpha + y$. We denote by $B_H(z,r)$ the Hilbert ball of center $z\in \R^n$ and radius $r\in \R_+$, i.e., $B_H(z,r) \defi \left\{ x \in  \R^n \colon \HNorm{x-z} \leq r \right\}$. We also endow $\trop$ with the standard order $\leq$, which is extended to vectors entrywise.

Another algebraic structure that we will use in this paper is the {\em completed min-plus semiring} $\rminb$, which is the set $\R\cup \{+\infty \}\cup \{-\infty \}$ equipped with $(a,b) \mapsto \min \{ a, b\}$ as addition and $(a,b) \mapsto a + b$ as multiplication (with the convention $(-\infty)+ (+\infty) = (+\infty)+ (-\infty) = (+\infty)$). The corresponding matrix product for matrices with entries in $\rminb$ will be denoted by $\odot'$. Given $A\in \trop^{m\times n}$, the operator $A^\sharp : \rminb^m \mapsto \rminb^n$ is defined by:
\[
A^\sharp (y) \defi (-A^\top) \odot' y \; ,
\]
where $A^\top$ denotes the transpose of $A$. The operator $A^\sharp$ will be called the {\em adjoint} of $A$, being an adjoint in a categorical sense
as it satisfies the following property:
\begin{equation}\label{ResiduationProp}
A \odot x \leq y \makebox{ if and only if } x \leq A^\sharp (y) \; ,
\end{equation}
for any $y\in \rminb^m$ and $x\in \trop^n$. 
 
\subsection{Stochastic zero-sum games with mean payoff}
\label{subsec:mpg}

A \emph{stochastic mean payoff game} can be specified by two matrices $A \in \trop^{m \times n}$ and $B \in \trop^{m \times q}$, and a row-stochastic matrix $P \in [0,1]^{q \times n}$, where $m,n,q\geq 1$. 
The rules of the game are as follows. Two players, called Max and Min, control disjoint sets of states, respectively indexed by $[m]$ and $[n]$, and alternatively move a pawn over these states. When the pawn is located on a state $j \in [n]$ of Player Min, she selects a state $i \in [m]$ such that $A_{ij} \neq \zero$, moves the pawn to state $i$ and pays to Player Max the amount $-A_{ij}$. When the pawn is on a state $i \in [m]$ of Player Max, he selects a state $k \in [q]$ such that $B_{ik} \neq -\infty$, moves the pawn to state $k$ and receives from Player Min the payment $B_{ik}$. Finally, at state $k$ the pawn is moved by {\em nature} to state $l \in [n]$ with probability $P_{kl}$. 

We shall make the following finiteness assumption, which assures that players Max and Min have at least one move available in each state.
\begin{assumption}\label{assum-fin}
Every row of $B$ has at least one finite entry, and the same is true for every column of $A$. 
\end{assumption}

A \emph{(positional) strategy} for Player Min is a function $\sigma \colon [n] \to [m]$ such that $A_{\sigma(j) j} \neq \zero$ for all $j$. Similarly, a \emph{(positional) strategy} for Player Max is a function $\tau \colon \oneto{m} \to \oneto{q}$ such that $B_{i \tau(i)} \neq \zero$ for all $i$. 
If Min and Max play according to the strategies $\sigma$ and $\tau$, and start from state $j_0 \in \oneto{n}$, the movement of the pawn is described by a Markov chain over the disjoint union $\oneto{m} \uplus \oneto{n} \uplus \oneto{q}$. Then, the \emph{payoff (of Player Max)} is defined as
 the average payoff
\[
g_{j_0}(\sigma, \tau) \defi \lim_{N \to +\infty} \mathbb{E}_{\sigma \tau} \Bigl(\frac{1}{N} \sum_{p = 1}^N (-A_{i_p j_p} + B_{i_p k_p}) \Bigr) \, ,
\]
where $\mathbb{E}_{\sigma \tau}$ refers to the expectation over the trajectories $j_0, i_0, k_0, j_1, i_1, k_1, \dots$, with respect to the probability measure
determined by these strategies. The objective of Players Min and Max is to find a strategy which respectively minimizes and maximizes the payoff. 
Liggett and Lippman~\cite{liggett_lippman} showed that there exists a pair of optimal strategies $(\sigma^*, \tau^*)$, which satisfies
\[
g_j(\sigma^*, \tau) \leq g_j(\sigma^*, \tau^*) \leq g_j(\sigma, \tau^*) \, ,
\]
for every initial state $j \in [n]$ and pair of Min/Max strategies $(\sigma, \tau)$. In this case, the quantity $g_j(\sigma^*, \tau^*)$ is referred to as the \emph{value} of the game when starting from state $j$. The state $j$ is said to be \emph{winning (for Player Max)} when the associated value is nonnegative. It is said to be \emph{strictly} winning for the same player if the associated value is positive. A dual terminology applies to Player Min.

With any such a game is associated a {\em Shapley operator}, which is the map $F:\trop^n\to \trop^n$ defined by
\begin{align}
F = A^\sharp \circ B \circ P \enspace, 
\label{e-niceshapley}
\end{align}
i.e., $F(x)=A^\sharp(B\odot (Px))$, where $Px$ denotes the usual
matrix-vector product of $P$ and $x$. The finiteness assumption
(\cref{assum-fin}) on the entries of the matrices $A,B$ imply that $F$ preserves both $\trop^n$ and $\R^n$. It is convenient %
to consider the vector $v^k \defi F^k(0)$, for $k\in \N$, where $F^k
= F\circ \dots\circ F$ denotes the $k$th iterate of $F$. The $j$th entry
$v^k_j$ represents the value of the game in finite horizon $k$ with initial
state $j$, associated with the same data. The vector
\[
\chi(F)\defi 
\lim_{k\to\infty} {v^k}/k 
= \lim_{k\to \infty} F^k(0)/k
\]
is known as the {\em escape rate vector} of $F$.
We shall recall in \cref{sec:perron_frobenius} why this escape rate exists.
It is known that 
\begin{align}
g_j(\sigma^*,\tau^*) = \chi_j(F)\enspace , 
\label{e-dynamic}
\end{align}
i.e., the value of the mean payoff game coincides with the limit
of the mean value per time unit of the finite horizon game, as
the horizon tends to infinity.
In this way, solving a mean payoff games reduces to
a dynamical systems issue: computing the 
escape rate vector of a Shapley operator.

Mean payoff games can be defined in different guises:
this only changes the explicit form of the Shapley
operator, without impact on the complexity of the problem,
as shown by the following remark.

\begin{remark}\label{rk-key}
Here, we assumed that Players Min, Max, and nature, play successively,
in this order.  Starting with Player Max, instead of Min, while keeping
the same circular order, would result
in replacing the Shapley operator $F$ by its cyclic conjugate
\begin{align}
F^{\mathrm{cyc}}(y) =  B\odot (P  (A^\sharp(y)))
\label{e-cyclic}
\end{align}
defined on $\trop^m$.
If the order was changed in a non cyclic way, nature playing for instance after
Min and before Max, then, the original Shapley operator $F$ would be replaced by:
\begin{align}
\bar{F}(x) = \bar{A}^\sharp (\bar{P} (\bar{B}\odot x)) \label{e-transp}
\end{align}
for some matrices $\bar A\in \trop^{q\times n}$, $\bar P\in [0,1]^{q\times m}$, and $\bar B\in \trop^{m\times n}$. It is also convenient to consider the effect of Players Max and Min swapping
their roles in the original game. This would amount to replacing $F$ by:
\begin{align}
\tilde{F}(x) = -F(-x) = A^\top \odot ((B^{\top})^{\sharp}( P x)) \; ,
\label{e-swap}
\end{align}
recalling that $\cdot^\top$ denotes the transposition.
Observe that $\chi(\tilde{F})=-\chi(F)$.
Moreover, if $F$ can be factored as $G\circ H$, $G$ and $H$ being any compositions of maps of the form $A^\sharp$, $y\mapsto B\odot y$, or $z\mapsto Pz$,
it can be shown
that $\chi(F)= \hat{G}(\chi(H\circ G))$ where $\hat{G}(x) \defi \lim_{s\to\infty}
s^{-1}G(sx)$ is the recession function of $G$, whose evaluation
is straightforward. Hence, one can recover the escape rate vector
of $G\circ H$ from the escape rate vector of its cyclic conjugate 
$H\circ G$, and vice versa.
Therefore, for an operator given in any of the forms~\cref{e-niceshapley,e-cyclic,e-transp,e-swap}, the complexity 
of computing the escape rate is independent of the choice of the special form.
\end{remark}

\subsection{Zero-sum games associated with nonarchimedean semidefinite programs}\label{subsec:correspondence}

The correspondence between semidefinite feasibility problems for spectrahedral cones and stochastic mean payoff games is given in the next theorem:
\begin{theorem}\label{th:correspondence}
With every spectrahedral cone $\bo \sC$ of the form~\cref{puiseux_spectra_cone} is associated a stochastic mean payoff game that satisfies the following property:
if the valuation of the entries of the matrices $\bo Q^{(i)}$ are chosen in a generic way, then $\bo \sC$ is nontrivial if and only if at least one state in the associated game is winning.
\end{theorem}

This correspondence is established in~\cite{issac2016jsc} by considering the following problem:
\begin{equation*}
\PP(F) \colon \text{does there exist} \; x \in \trop^n \; \text{such that } x \neq \mpzero 
 \; \text{and} \; x \leq F(x)?
\end{equation*}
where $F \colon \trop^n \to \trop^n$ is the Shapley operator of the game associated with the spectrahedral cone $\sC$. 
This problem is said to be \emph{feasible} when it admits a solution, and \emph{infeasible} otherwise. We point out that $\PP(F)$ is feasible if, and only if, the associated stochastic mean payoff game has a winning state. Equivalently, this amounts to the fact that the set 
\begin{equation}\label{FeasibleSet}
\sS(F) \defi \{ x\in \trop^n \colon x \leq F(x)\} 
\end{equation}
is nontrivial, meaning that it is not reduced to the point $\mpzero$. 

The correspondence between nonarchimedean semidefinite
programming and stochastic games is simpler to present
if we assume that the matrices $\bo Q^{(1)}, \dots, \bo Q^{(n)}$ are \emph{(negated) Metzler matrices}, which means that their off-diagonal entries are nonpositive. In this case, if the genericity assumption of \cref{th:correspondence} is satisfied, $\sS(F)$ is precisely the image under the valuation map of $\sC$.
Similarly, we can consider the problem: 
\begin{equation*}
\PP_\R(F) \colon \text{does there exist} \; x \in \R^n \; \text{such that} \; x \strictleq F(x)?
\end{equation*}
where $y \strictleq z$ stands for the fact that $y_i < z_i$ for all $i$. This problem is feasible if, and only if, the set $\sS(F)$ is \emph{strictly nontrivial}, meaning that there exists $x \in \R^n$ such that $x \strictleq F(x)$. This corresponds to the property where every state of the game has a positive value. 

The Shapley operator and the associated feasibility problems $\PP(F)$ and $\PP_\R(F)$ provide further conditions under which game algorithms are directly applicable to solve nonarchimedean feasibility problems, disregarding the genericity conditions of \cref{th:correspondence}. 
\begin{theorem}\label{th:sufficient_conditions}
For any Metzler matrices $\bo Q^{(1)}, \dots, \bo Q^{(n)}$, we have: 
\begin{enumerate}[(i)]
\item if $\PP(F)$ is infeasible, or equivalently, $\sS(F)$ is trivial, then $\bo \sC$ is trivial.
\item if $\PP_\R(F)$ is feasible, or equivalently, $\sS(F)$ is strictly nontrivial, then $\bo \sC$ is strictly nontrivial, meaning that there exists $\bo x \in \puiseux_{>0}^n$ such that the matrix $\bo x_1 \bo Q^{(1)} + \dots + \bo x_n \bo Q^{(n)}$ is positive definite.
\end{enumerate}
\end{theorem}
Following the analogy with the classical condition number in linear programming (see, e.g., \cite{Renegar1995}), we are interested in finding a numerical quantity measuring (the inverse of) the distance to triviality when the instance is nontrivial or to nontriviality when it is trivial. In more details, we define the {\em condition number} $\cond(F)$ of the above problem $\PP(F)$ by:
\begin{equation}\label{DefCondNumber1}
(\inf \{\|u\|_\infty \colon u \in \R^n \, , \; \PP(u+F) \; \text{is infeasible} \} )^{-1}  
\end{equation}
if $\PP(F)$ is feasible, and 
\begin{equation}\label{DefCondNumber2}
(\inf \{\|u\|_\infty \colon u \in \R^n \, , \; \PP(u+F) \; \text{is feasible} \} )^{-1}   
\end{equation}
if $\PP(F)$ is infeasible (with the convention $0^{-1} = +\infty$). Here, $u + F$ stands for the map $x \mapsto u + F(x)$, where the addition is understood entrywise, and $\| \cdot \|_\infty$ stands for the sup-norm, \ie, $\| u \|_\infty \defi \max_i |u_i|$. The condition number $\cond_{\R}(F)$ of the problem $\PP_\R(F)$ is defined in the same way as in~\cref{DefCondNumber1} and~\cref{DefCondNumber2} but replacing $\PP$ by $\PP_\R$. 

\begin{remark}
Looking for additive perturbations of the form $u+F$ is a canonical approach; such perturbations have been already used to reveal the ergodicity properties of the game~\cite{ergodicity_conditions}. This is also the finite dimensional analogue of perturbing the Hamiltonian of a Hamilton--Jacobi PDE by adding a potential~\cite{rifford}.
\end{remark}
\begin{remark}
In~\cite{issac2016jsc}, the Shapley operator associated with a nonarchimedean
SDP feasibility problem is written as $A^\sharp \circ P \circ B$ instead of $A^\sharp \circ B \circ P$. 
As there are reductions between the games corresponding to both forms
(as discussed in~\cref{rk-key}),
we consider here a Shapley operator in the latter form.
This is more suitable to state the complexity estimates in~\cref{sec-complexity}. 
\end{remark}

\section{Preliminary results of nonlinear Perron--Frobenius theory}

\label{sec:perron_frobenius}

In this section, we recall some elements of nonlinear Perron--Frobenius theory which will be used to study the condition numbers introduced above. To do so,
we next axiomatize essential properties of the Shapley operators considered in~\cref{subsec:mpg}, following the ``operator approach'' of stochastic games~\cite{RS01,Ney03}.

A self-map $F$ of $\trop^n$ is said to be {\em order-preserving} when
\[
x\leq y \implies F(x) \leq F(y) \; \text{for all} \; x,y\in \trop^n \, ,
\]
and {\em additively homogeneous} when 
\[ 
F(\lambda  + x) =\lambda + F(x) \; \text{for all} \; \lambda \in \trop \; \text{and} \; x\in \trop^n \, .
\]
We point out that any order-preserving and additively homogeneous self-map $F$ of $\trop^n$ that preserves $\R^n$ is nonexpansive in the sup-norm, meaning that
\[
\|F(x) - F(y)\|_\infty \leq \|x - y\|_\infty \; \text{for all} \; x,y\in \R^n \; .
\]

Given an order-preserving and additively homogeneous self-map $F$ of $\trop^n$, the vectors $x\in \trop^n$ satisfying $x\leq F(x)$ can be thought of as the nonlinear analogues of subharmonic functions. 
A central role in determining the existence of such vectors is played by the limit 
$\chi(F)=\lim_{k \rightarrow \infty} (F^k(x)/k)$, for $x\in \R^n$. 
When this limit exists, it can be shown to be independent of the choice of $x\in \R^n$, and so it coincides with the escape rate vector $\chi (F)$ of $F$.  
The following theorem of Kohlberg implies that the limit does exist
 when $F$ preserves $\R^n$ and its restriction to $\R^n$ is {\em piecewise affine} (meaning that $\R^n$ can be covered by finitely many polyhedra such that $F$ restricted to any of them is affine). 

\begin{theorem}\cite{kohlberg}\label{Theo:Kohlberg}
A piecewise affine self-map $F$ of $\R^n$ that is nonexpansive in any norm admits an invariant half-line, meaning that there exist $z, w\in \R^n$ such that 
\[
F(z+\beta w) = z + (\beta +1) w
\]  
for any $\beta \in \R$ large enough. In particular, the escape rate vector $\chi(F)$ exists, and is given by the vector $w$.
\end{theorem}
Kohlberg's theorem applies to Shapley operators of stochastic mean payoff games with finite state and action spaces and perfect information. 
Indeed, the Shapley operator~\cref{e-niceshapley} of the game described in \cref{subsec:mpg} is order-preserving and additively homogeneous, and its restriction to $\R^n$ is piecewise affine. 

For a general order-preserving and additively homogeneous self-map of $\trop^n$,
the escape rate vector may not exist. We can still, however, recover information about the sequences $(F^k(x)/k)_k$ through the Collatz--Wielandt numbers of $F$. Assuming that $F$ is a continuous, order-preserving, and additively homogeneous self-map $F$ of $\trop^n$, we define the {\em upper Collatz--Wielandt number} of $F$ by:
\begin{equation}\label{DefUCW}
\ucw (F) \defi \inf \{ \mu \in \R \colon \exists z \in \R^n , F(z)\leq \mu + z \} \, ,
\end{equation}
and the {\em lower Collatz--Wielandt number} of $F$ by:
\begin{equation}\label{DefLCW}
\lcw (F) \defi \sup \{ \mu \in \R \colon \exists z \in \R^n , F(z)\geq \mu + z \} \, .
\end{equation}
A relation between the escape rate vector and the upper Collatz--Wielandt number is given in the next theorem, which is derived
in~\cite{polyhedra_equiv_mean_payoff} from a minimax result of Nussbaum~\cite{nussbaum}. 

\begin{theorem}\cite[Lemma~2.8]{polyhedra_equiv_mean_payoff}\label{Theo:UCW}
Let $F$ be a continuous, order-preserving, and additively homogeneous self-map of $\trop^n$. Then, 
\begin{align*}
 \lim_{k \rightarrow \infty} \mytop (F^k(x)/k) = \ucw (F) =
 \sup \{ \mu \in \trop \colon \exists z \in \trop^n , \; z\neq \mpzero, \; F(z)\geq \mu + z \}
\end{align*}
for any $x\in \R^n$.
\end{theorem}

It is known that an order-preserving and additively homogeneous self-map of $\R^n$ admits a unique continuous extension to $\trop^n$, see~\cite{burbanks_nussbaum_sparrow}. Then, as noted in~\cite[Remark~2.10]{polyhedra_equiv_mean_payoff}, the previous result can be dualized when $F$ preserves $\R^n$. 

\begin{corollary}\label{Coro:LCW}
Let $F$ be a continuous, order-preserving, and additively homogeneous self-map of $\trop^n$ that preserves $\R^n$. Then, 
\[
\lim_{k \rightarrow \infty} \mybot (F^k(x)/k) = \lcw (F)
\]
for any $x\in \R^n$. 
\end{corollary}
As a consequence, when the escape rate vector exists, we simply have
\[
\ucw(F) = \mytop(\chi(F)) \quad \text{and} \quad \lcw(F) = \mybot(\chi(F)) \, .
\]
Specializing this to the case where $F$ is the Shapley operator of a game, the quantities $\ucw(F)$ and $\lcw(F)$ respectively correspond to the greatest and smallest values of the states for the mean payoff problem.

In the sequel, we will consider especially the situation in which there is a vector $v \in \R^n$ and a scalar $\lambda\in\R$ such that
\begin{align}
F(v) = \lambda + v  \, .
\label{e-ergodic}
\end{align}
The scalar $\lambda$, which is
unique, is known as the {\em ergodic constant},
and~\cref{e-ergodic} is referred to as the {\em ergodic equation}.
We will denote this scalar by $\rho(F)$ as it is a nonlinear extension of the 
spectral radius. 
The vector $v$
is known as a {\em bias},
or a {\em potential}. 
It is easily seen that if $F$ admits such
a bias vector, then
\[
\lcw(F)=\ucw(F) = \rho(F) \, ,
\]
and the condition that $\rho(F) \geq 0$
means that the game is winning for every initial state.
The existence of a bias vector is guaranteed by certain ``ergodicity''
assumptions~\cite{ergodicity_conditions}. 

\section{Metric geometry properties of condition numbers}\label{sec-metric}
\subsection{Condition numbers vs Collatz--Wielandt numbers, and duality}\label{sec-duality}

We point out that the definitions given in \cref{subsec:correspondence} of the condition numbers $\cond(F)$ and $\cond_\R(F)$ can be generalized to any continuous, order-preserving, and additively homogeneous self-map $F$ of $\trop^n$. The next proposition provides a characterization of these condition numbers in terms of the Collatz--Wielandt numbers of $F$. 
\begin{proposition}\label{Lem:CondNumCWRel}
Let $F$ be a continuous, order-preserving, and additively homogeneous self-map of $\trop^n$. Then, 
\[
\cond_\R(F)= |\lcw(F)|^{-1}
\; \text{and} \; 
\cond(F)= |\ucw(F)|^{-1}.
\]
\end{proposition} 

We define the dual of the mean payoff game of \cref{subsec:correspondence} as the one whose Shapley operator is $F^\ast = (B^\top)^\sharp \circ P \circ A^\top$. The following theorem will allow us to relate $\PP_\R (F)$ with $\PP_\R (F^\ast)$ and $\PP (F^*)$. 

\begin{theorem}[Duality theorem]\label{Th:RelCWNumbers}
Let $F=A^\sharp \circ B \circ P$ and $F^\ast = (B^\top)^\sharp \circ P \circ A^\top$, where $A\in \trop^{m\times n}$ and $B\in \trop^{m\times q}$ satisfy \cref{assum-fin}, $A$ has at least one
finite entry 
per row, and $P\in \R^{q\times n}$ is a row-stochastic matrix. Then, 
\[
\ucw(F^\ast) = - \lcw(F) \; .
\] 
\end{theorem}

As a consequence of \cref{Th:RelCWNumbers}, we obtain: 

\begin{corollary}\label{Cor:Dual}
Let $F=A^\sharp \circ B \circ P$ and $F^\ast = (B^\top)^\sharp \circ P \circ A^\top$, where $A\in \trop^{m\times n}$ and $B\in \trop^{m\times q}$ satisfy \cref{assum-fin}, $A$ has at least one
finite entry 
per row, and $P\in \R^{q\times n}$ is a row-stochastic matrix. Then, 
\begin{enumerate}[(i)]

\item\label{CorDualItem:CondNumEqualDual}
The condition number of $\PP_\R(F)$ coincides with the condition number of $\PP (F^\ast)$.

\item\label{CorDualItem:StrictFeasVSFeas}
Either $\PP (F^*)$ is feasible or $\PP_\R(F)$ is feasible. 

\item\label{CorDualItem:StrictlyFeasible}
Only one of the problems $\PP_\R (F)$ and $\PP_\R (F^\ast)$ can be feasible. 

\end{enumerate}
\end{corollary}

\subsection{A geometric characterization of condition numbers}

In this section, we study the {\em inner radius} of the feasible sets of games, that is, given the Shapley operator $F\colon \trop^n \to \trop^n$ of a game, we study the maximal radius of a Hilbert ball contained in the set~\cref{FeasibleSet}.

We start with the following simple lemma. 

\begin{lemma}\label{LemmaSufBoundRadius}
Let $F$ be an order-preserving and additively homogeneous self-map of $\trop^n$. Assume $z\in \R^n$ and $r\in \R_+$ are such that $r \leq \mybot(F(z)-z)$. Then, the Hilbert ball $B_H(z,r)$ is contained in $\sS(F)$. 
\end{lemma}

For the condition in the previous lemma to be also necessary for the inclusion to hold, we need an additional assumption on $F$. 

\begin{definition}
An order-preserving and additively homogeneous self-map $F$ of $\trop^n$ is said to be {\em diagonal free} when $F_i(x)$ is independent of $x_i$ for all $i\in \oneto{n}$. In other words, $F$ is diagonal free
if for all $i\in\oneto{n}$, and for all $x,y\in\R^n$ such that $x_j=y_j$
for $j\neq i$, we have $F_i(x)=F_i(y)$. 
\end{definition}

\begin{lemma}\label{LemmaNecBoundRadius}
When $F$ is diagonal free, for any $z\in \R^n$ and $r\in \R_+$ the Hilbert ball $B_H(z,r)$ is contained in $\sS(F)$ only if $r \leq \mybot(F(z)-z)$.
\end{lemma}

If $F$ is not diagonal free, the conclusion of \cref{LemmaNecBoundRadius} does not necessarily hold, as shown in the next example.

\begin{example}
Let us consider the order-preserving and additively homogeneous map $F=A^\sharp \circ B$, where $A=\begin{pmatrix} 0 & 0 & 0 \end{pmatrix}$ and $B=\begin{pmatrix} -1 & 0 & -1 \end{pmatrix}$. 
Then, for $z=\begin{pmatrix} 0 & 3 & 0 \end{pmatrix}^\top$, it can be verified that $B_H(z,3)\subset \sS(F) = \left\{ x \in \R^3 \colon x \leq A^\sharp \circ B (x) \right\} = \left\{ x \in \R^3 \colon A\odot x \leq B\odot x \right\} $. However, we have 
\[
F(x)=\begin{pmatrix} 
\max \{x_1-1,x_2,x_3-1\} \\ 
\max \{x_1-1,x_2,x_3-1\} \\ 
\max \{x_1-1,x_2,x_3-1\} \end{pmatrix}\; ,
\] 
and so $\mybot(F(z)-z)=0$. 
\end{example}

As a consequence of \cref{LemmaSufBoundRadius,LemmaNecBoundRadius}, we obtain: 

\begin{theorem}\label{ThMaximalRadius}
Let $F$ be a diagonal free self-map of $\trop^n$. Then, $\sS(F)$ contains a Hilbert ball of positive radius if and only if $\lcw(F)> 0$. Moreover, when $\sS(F)$ contains a Hilbert ball of positive radius, the supremum of the radii of the Hilbert balls contained in $\sS(F)$ coincides with $\lcw(F)$.
\end{theorem}

Sergeev established in~\cite{Sergeev2007} a characterization of the inner radius of {\em polytropes}, which corresponds to the special case of \cref{ThMaximalRadius}
in which $F$ is the Shapley operator of a game with only one player and deterministic transitions.

\begin{remark}
The condition in \cref{ThMaximalRadius} is not too restrictive. Indeed, it can be shown that in most cases of interest, if the Shapley operator $F$ of a mean payoff game is not diagonal free, one can construct another mean payoff game such that its Shapley operator is diagonal free and the inner radius of its feasible set coincides with the one of $\sS(F)$. 
\end{remark}

\section{Bounding the complexity of value iteration by the condition numbers}
\label{sec-complexity}
In this section, $F$ is an order-preserving and additively homogeneous self-map of $\trop^n$ which preserves $\R^n$. 
We also assume that $F$ admits a bias vector $v\in\R^n$, as in~\cref{e-ergodic}.

\subsection{A universal complexity bound for value iteration}
The most straightforward idea to solve
a mean payoff game is probably value iteration:
we infer whether or not the mean payoff
game is winning by solving the finite horizon game,
for a large enough horizon. This is formalized in \cref{algo-vi}.

\begin{figure}[t]
\begin{small}
\begin{algorithmic}[1]
\Procedure {ValueIteration}{$\shapley$}
\\ \Lcomment{$\shapley$ a Shapley operator from $\R^n$ to $\R^n$ }
\\\Lcomment{The algorithm will report whether Player Max or Player Min wins the mean payoff game represented by $F$}
\State $u \coloneqq 0 \in \R^n$ 
\While{$\mytop (u)>0$ and $\mybot (u)<0$}
\label{step-while}
\label{state-iter} $u \coloneqq \shapley(u)$  
\Lcomment{At iteration $\ell$, $u=F^\ell(0)$ is the value vector of the game in finite horizon $\ell$}
\EndWhile
\If{$\mytop (u)<0$}\label{step-termmax}
\Return ``Player Min wins''
\Else \label{step-termmin}
\ \Return ``Player Max wins''
\EndIf
\EndProcedure
\end{algorithmic}
\end{small}
\caption{Basic value iteration algorithm.} \label{algo-vi}
\end{figure}

We next show, in \cref{Th:NitsBound}, that this value iteration
algorithm terminates and is correct, provided
the mean payoff of the game is nonzero (i.e., $\rho(F) \neq 0$), and the
operations are performed in exact arithmetic.
We shall see in \cref{cor-perturb,cor-approx} that these two restrictions
can be eliminated, at the price of an increase 
of the complexity bound.

It is convenient to introduce the following metric estimate, which
represents the minimal Hilbert's seminorm of a bias vector
\[
R(F) \defi \inf \left\{ \HNorm{u} \colon u \in \R^n , \; F(u)=\rho (F) + u \right\} \, .  
\]
Since $F$ is assumed to have a bias vector $v\in \R^n$, 
we have $R(F)\leq \|v\|_H<\infty$ and
$\rho(F)=\lcw(F)=\ucw(F)$. Hence, by \cref{Lem:CondNumCWRel}, 
\begin{equation*}
|\rho(F)|^{-1}=|\lcw(F)|^{-1} = |\ucw(F)|^{-1} = \cond_\R(F)= 
\cond(F)
\, .\end{equation*}
We shall denote by $\cond(F)$ this common quantity.

Note that $|\rho(F)|$ has a remarkable interpretation,
as the value of an auxiliary game, in which there is an initial stage,
at which Player Max can decide either to keep his role or to swap
it with the role of Player Min.
Then, the two players play the mean payoff game as usual. 
Swapping roles amounts to 
replacing $F$ by the Shapley operator $\tilde{F}(x) \defi -F(-x)$.
Observe also that $\rho(\tilde{F})=-\rho(F)$ as noted in \cref{rk-key}.
Hence, the value of this modified game is precisely $\max(\rho(F),\rho(\tilde{F})) = |\rho(F)|$. 

The following result bounds the complexity of value iteration
in terms of $R(F)$ and of the condition number $\cond(F)$.
\begin{theorem}\label{Th:NitsBound}
Suppose that the Shapley operator $F$ has a bias vector and that the ergodic constant $\rho(F)$ is nonzero. Then, procedure
\textsc{ValueIteration} terminates after
\[ 
N_{\mathrm{vi}} \leq R(F)\cond(F) \;
\]
iterations and returns the correct answer. 
\end{theorem}

\subsection{Bounding the condition number and the bias vector of a stochastic mean payoff game}

We next bound the condition number $|\rho(F)|^{-1}$, and the metric estimate
$R(F)$, when $F$ is a Shapley operator of a stochastic game
with perfect information and finite action spaces.
As in \cref{subsec:mpg},
we assume that 
\begin{align}
\shapley =A^\sharp \circ B \circ P
\label{e-def-shapley}
\end{align}
where $A\in \trop^{\maxstates\times \minstates}$ has at least one finite entry 
per column, $B\in \trop^{\maxstates\times \chancestates}$ has at least one finite entry per row, and $P\in \R^{\chancestates\times \minstates}$ is a row-stochastic matrix. 
To obtain explicit bounds, we will assume that the finite entries of the matrices $A$ and $B$ are {\em integers},
and we set 
\begin{equation}\label{e-defw}
\nonumber
W\defi \max \left\{| A_{ij} - B_{ih} | \colon A_{ij} \neq \mpzero, \, B_{ih} \neq \mpzero, \, i \in \oneto{m}, \, j\in \oneto{n}, \, h \in \oneto{q} \right\}\, .
\end{equation}
This is not more special than assuming that the finite entries of $A$ and $B$
are rational numbers (we may always rescale rational payments so that they become integers).
We also assume that the probabilities $\transition_{\state \stateII}$ are rational, and that they
have a common denominator $\comd \in \N_{> 0}$, $\transition_{\state \stateII} = \potential_{\state \stateII}/\comd$, where $\potential_{\state \stateII} \in [\comd]$ for all $\state\in [\chancestates]$ and $\stateII \in [\minstates]$.

We say that a state $i\in[\chancestates]$ is {\em nondeterministic} 
if there are at least two indices $l,l'\in[\minstates]$ such that
$P_{il}>0$ and $P_{i l'}>0$.

The following lemma improves an estimate in~\cite{boros_gurvich_makino}. 
\begin{lemma}\label{est_invariant_measure}
Suppose that a Markov chain with $\nstates$ states is irreducible, and that the transition probabilities are rational numbers whose denominators divide an integer $M$. Let $\rstates\leq \minstates$ denote the number of states
with at least $2$ possible successors. Let $\invmes \in (0,1]^{\minstates \times \minstates}$ denote the invariant measure of the chain. Then, the least common denominator of the rational numbers $(\invmes_{\state})_{\state \in \states}$ is not greater than $\nstates\comd^{\min\{\rstates,\nstates-1\}}$.
\end{lemma}
We deduce the following result. 
\begin{theorem}\label{th-cond}
Let $F$ be a Shapley operator as above, still supposing that $F$ has a bias vector and that $\rho(F)$ is nonzero. If $\rstates$ is the number of nondeterministic states of the game, then $\cond(F) \leq \nstates\comd^{\min\{\rstates,\nstates-1\}}$.
\end{theorem}

To bound $R(F)$, we use the following idea.
For $0<\alpha<1$, let $v_\alpha$ denote the value
of the discounted game associated with $F$, meaning
that $v_\alpha = F(\alpha v_\alpha)$. 
Since $F$ 
represents a zero-sum game with perfect information
and finite state and action spaces, it is known that
$v_\alpha$ has a Laurent series expansion in powers
of $(1-\alpha)$ with a pole of order at most~$1$
at $\alpha=1$, see~\cite{kohlberg}. 
We can deduce from this that the limit of $v_\alpha - \rho(F) /(1-\alpha)$
as $\alpha \to 1^-$ exists and that it is a bias,
which we call the {\em Blackwell bias}. By working out the limit,
we arrive at the following estimate.
\begin{theorem}\label{th-bias}
Let $\shapley$ be the Shapley operator in~\cref{e-def-shapley},
still supposing that it has a bias vector,
and let $\bias^*$ be its Blackwell bias.
Then,
\[ R(\shapley)\leq \hilb{\bias^*} \le 10\nstates^{2} \Payoff \comd^{\min\{\rstates, \nstates - 1\}} \, .\]
\end{theorem}
By combining \cref{th-bias,th-cond}, 
we arrive at the following.
\begin{corollary}\label{coro-totalbound}
Let $\shapley$ be the Shapley operator in~\cref{e-def-shapley},
still supposing that it has a bias vector 
and that $\rho(F)$ is nonzero. Then,
procedure
\textsc{ValueIteration} stops
after
\begin{align}
N_{\mathrm{vi}}\leq 10\nstates^{3} \Payoff \comd^{2\min\{\rstates, \nstates - 1\}} 
\label{e-firstbound}
\end{align}
iterations and correctly decides which of the two players is winning.
\end{corollary}

We next show that when specialized
to deterministic games, 
the universal estimate of \cref{Th:NitsBound}
gives precisely the complexity bound of Zwick--Paterson~\cite{zwick_paterson}.

\begin{lemma}\label{Lem:BoundRDet}
Let $F=A^\sharp \circ B$, where $A,B\in \trop^{m\times n}$, and suppose that there exists $v\in \R^n$ such that $F(v)= \rho(F) +v$. Then
\[
R(F) \leq (n-1) (|\rho (F)| + W) \, ,
\]
where $W$ is defined as in~\cref{e-defw}, setting $q=n$.
\end{lemma}
For deterministic games with integer payments, the mean payoff
is given by the average weight of a circuit,
which has length at most $\minstates$. It follows that $|\rho(F)|\geq 1/n$, unless $\rho(F)=0$. Note also that $\rho(F)\leq W$. By applying \cref{th-cond}, we arrive at the following bound for the number of iterations $N_{\mathrm{vi}}$ of the algorithm in \cref{algo-vi}.

\begin{corollary}[Compare with~\cite{zwick_paterson}]
 Let $F=A^\sharp \circ B$ be the Shapley operator
of a deterministic game, where the finite entries of $A,B\in \trop^{m\times n}$ are integers. If there exists $v\in \R^n$ such that $F(v)= \rho(F)+v$ with $\rho(F)\neq 0$,
then
\[ 
N_{\mathrm{vi}} \leq 2 n^2 W \; .
\]
\end{corollary}

The assumption $\rho(F)\neq 0$ that is used in \cref{Th:NitsBound}
can be relaxed, by appealing to the following 
perturbation and scaling argument. This leads
to a bound in which the exponents of $M$ and of $n$ are increased. 
\begin{corollary}\label{cor-perturb}
Let $\mu \defi\nstates\comd^{\min\{\rstates,\nstates-1\}}$.
Then, procedure
\textsc{ValueIteration}, applied
to the perturbed and rescaled Shapley operator $1+2\mu F$, 
satisfies
\[ 
N_{\mathrm{vi}} \leq 21\nstates^{4} \Payoff \comd^{3\min\{\rstates, \nstates - 1\}} 
\] 
iterations, and this holds unconditionally.
If the algorithm reports that Max wins, then 
Max is winning in the original mean payoff game.
If the algorithm reports that Min wins, then Min is strictly winning in the original mean payoff game.
\end{corollary}

The algorithm in \cref{algo-vi} can be adapted to work
in finite precision arithmetic. Consider the variant
of the main body of this algorithm, given in \cref{algo-vi2}. 
We assume that each evaluation of the Shapley operator $F$
is performed with an error of at most $\epsilon>0$ in the sup-norm.
\begin{figure}[t]
\begin{small}
\begin{algorithmic}
\State $u \coloneqq 0 \in \R^n$, $\ell\coloneqq 0\in \N$, $\epsilon\in \R_{>0}$ 
\While{$\ell\epsilon+\mytop (u)\geq 0$ and $-\ell\epsilon+\mybot (u)\leq 0$}
\label{step-whilenew}
\State 
\label{state-iternew} $u \coloneqq \shapley(u)$; $\ell\coloneqq \ell+1$  
\Lcomment{The operator $\shapley$ is evaluated in approximate arithmetic, so that $F(u)$ is at most at distance $\epsilon$ in the sup-norm from its true value.}
\EndWhile
\If{$\ell\epsilon+ \mytop (u)\leq 0$}\label{step-termmax2new}
\Return ``Player Min wins''
\EndIf
\If{$-\ell\epsilon+ \mybot (u)\geq 0$}\label{step-termmax3new}
\Return ``Player Max wins''
\EndIf
\end{algorithmic}
\end{small}
\caption{Modification of the basic value iteration algorithm to work in finite precision arithmetic.} 
\label{algo-vi2}
\end{figure}

\begin{corollary}\label{cor-approx}
Let $\shapley$ be the Shapley operator in~\cref{e-def-shapley},
still supposing that it has a bias vector
and that $\rho(F)$ is nonzero. 
Let $\mu \defi \nstates\comd^{\min\{\rstates,\nstates-1\}}$.
Then, for any $0<\epsilon\leq \mu^{-1}/3$, 
value iteration performed with a numerical precision of
$\epsilon$ at each step (i.e., the algorithm in \cref{algo-vi2}) stops
after
\begin{align}
 N_{\mathrm{vi}}\leq 30\nstates^{3} \Payoff \comd^{2\min\{\rstates, \nstates - 1\}} 
\label{e-boundapprox}
\end{align}
iterations and correctly decides which of the two players is winning.
\end{corollary}
Observe that~\cref{e-boundapprox} is the bound~\cref{e-firstbound}
multiplied by $3$. 

\section{Concluding remarks}
We introduced a notion of condition number for stochastic mean payoff games,
and bounded the complexity of value iteration
in terms of this condition number. 
Whereas condition numbers are familiar for 
problems over archimedean fields, this
leads to an appropriate notion of condition
number for {\em nonarchimedean} semidefinite
programming. In particular, our present results
explain, at least in part,
the perhaps surprising benchmarks of~\cite{issac2016jsc}, 
revealing that random nonarchimedean semidefinite feasibility 
instances with generic valuations can be simpler
to solve than their archimedean analogues. In some sense,
``good conditioning'' provides a quantitative version
of ``genericity,'' and most instances in~\cite{issac2016jsc}
are well conditioned. This raises the issue of evaluating the condition
number on random instances.
It is also an interesting question to investigate
whether the solution of nonarchimedean SDP could be used, in general,
to solve archimedean SDP, and vice versa.

\section*{Acknowledgement}
The second author thanks Vladimir Gurvich for enlightening discussions on the pumping algorithm of~\cite{gurvich} and its extension to stochastic games in~\cite{boros_gurvich_makino}.
\bibliographystyle{alpha}

\end{document}